\def\QTR#1#2{{\csname#1\endcsname #2}}

\documentclass{amsart}
\usepackage{amsfonts}
\usepackage{amscd}
\usepackage{thmdefs}

\theoremstyle{definition}
\theoremstyle{remark}
\numberwithin{equation}{section}

\begin{document}
\ifx\hyperref\undefined\else\errmessage{hyper disabled}\fi
\begin{flushright}
UAHEP 956\\
November 1995\\
\end{flushright}
\bigskip
\bigskip

\title{Two-parameter deformation of the Poincar\'e algebra}
\author{A. Stern}
\address{Dept. Physics and Astronomy, Univ. of Alabama, Tuscaloosa, AL 35487, USA.}
\email{astern@@ua1vm.ua.edu and iyakush3@@ua1vm.ua.edu}
\author{I. Yakushin}
\maketitle

\begin{abstract}
We examine a two-parameter ($\hbar ,$ $\lambda $) deformation of the
Poincar\`e algebra which is covariant under the action of $SL_q(2,C).$ When $%
\lambda \rightarrow 0$ it yields the Poincar\`e algebra, while in the $\hbar
\rightarrow 0$ limit we recover the classical quadratic algebra discussed
previously in \cite{ssy95}, \cite{sy95}. The analogues of the Pauli-Lubanski
vector $w$ and Casimirs $p^2$ and $w^2$ are found and a set of mutually
commuting operators is constructed.
\end{abstract}

\section{Introduction\ }

In \cite{ssy95}, \cite{sy95} we proposed a two-parameter deformation of
Poincar\'e algebra which transforms covariantly under the action of $%
SL_q(2,C)$. The algebra appears to be distinct from systems discussed
previously, e.g. in \cite{oswz92}, \cite{PSW93}, \cite{swz91}, \cite{sm1194}%
, \cite{sm895}, \cite{nl91}, \cite{rt94}, \cite{lr94}, \cite{bhos95} and it
has the advantage that it can be expressed very compactly in terms of two 2$%
\times $2 matrices which we denote by $P$ and $\Gamma .$ The corresponding
matrix elements $P_{ij}$ and $\Gamma _{ij}$ are operators analogous to
momentum and angular momentum, respectively. In order to have a ten
dimensional algebra, we imposed the following constraints on $P$ and $\Gamma
:$ $P$ is hermitian (with respect to some anti-involution operation $\dagger 
$ ), while $\Gamma $ satisfies a certain deformed unimodularity condition
which we define shortly. 

We shall examine this deformed Poincar\'e algebra in more detail in this
article. In particular, we shall be interested in obtaining its Casimir
operators, along with a complete set of commuting operators. Representations
of this algebra will be constructed in forthcoming paper.

The algebra of \cite{sy95} is given by: 
\begin{eqnarray}
\underset{12}{R}\underset{1}{P}\underset{12}{R}^{-1}\underset{2}{P} &=&%
\underset{2}{P}\underset{21}{R}^{-1}\underset{1}{P}\underset{21}{R},
\label{pp} \\
\underset{21}{R}^{-1}\underset{1}{\Gamma }\underset{21}{R}\underset{2}{%
\Gamma } &=&\underset{2}{\Gamma }\underset{12}{R}\underset{1}{\Gamma }%
\underset{12}{R}^{-1},  \label{gg} \\
\underset{12}{R}\underset{1}{\Gamma }\underset{12}{R}^{-1}\underset{2}{%
\overline{\Gamma }} &=&\underset{2}{\overline{\Gamma }}\underset{12}{R}%
\underset{1}{\Gamma }\underset{12}{R}^{-1},  \label{ggb } \\
\underset{21}{R}^{-1}\underset{1}{P}\underset{21}{R}\underset{2}{\Gamma } &=&%
\underset{2}{\Gamma }\underset{21}{R}^{-1}\underset{1}{P}\underset{12}{R}%
^{-1},  \label{pg}
\end{eqnarray}
where $\overline{\Gamma }=\Gamma ^{\dagger -1}$ . We use tensor product
notation, labels 1 and 2 denote different vector spaces, $\underset{1}{P}%
=P\otimes \Bbb{I},$ $\underset{2}{P}=\Bbb{I}\otimes P$ , etc. where $\Bbb{I}$
is 2$\times $2 unit matrix. The $R-$ matrix is given by 
\begin{equation}
\underset{12}{R}=q^{-1/2}\left( 
\begin{array}{cccc}
q & 0 & 0 & 0 \\ 
0 & 1 & 0 & 0 \\ 
0 & q-q^{-1} & 1 & 0 \\ 
0 & 0 & 0 & q
\end{array}
\right)  \label{R}
\end{equation}
and satisfies the quantum Yang-Baxter equation. Additional commutational
relations are obtained by hermitian conjugation. Here $q$ is a real number.

The $P-P$ commutational relations (\ref{pp}) are already known, see for
instance \cite{akr94}, \cite{akr93}. They are consistent with $P\ $being
hermitian. The commutational relations (\ref{gg}-\ref{pg}) are new . We
supplement them by imposing a deformed unimodularity condition on $\Gamma $
which we now explain. By expanding (\ref{pp}-\ref{pg}) in terms of matrix
elements it can be shown that the following quadratic combinations of $%
\Gamma ^{\prime }s$ and $\overline{\Gamma }^{\prime }s$ are in the center of
the algebra and we can therefore make them equal to 1: 
\begin{eqnarray}
det_{\frac 1q}(\Gamma ^T) &=&\Gamma _{11}\Gamma _{22}-q^2\Gamma _{21}\Gamma
_{12}=1,  \label{c3} \\
det_q(\overline{\Gamma }) &=&\overline{\Gamma }_{11}\overline{\Gamma }%
_{22}-\frac 1{q^2}\overline{\Gamma }_{12}\overline{\Gamma }_{21}=1.
\label{c4}
\end{eqnarray}
(\ref{c3}-\ref{c4}) are not independent: one can be obtained from the other
by applying hermitian conjugation and using the expressions for $\Gamma
^{-1} $ and $\overline{\Gamma }^{-1}$%
\begin{eqnarray}
\overline{\Gamma }^{-1} &=&\frac{\left( 
\begin{array}{cc}
\overline{\Gamma }_{22} & -\frac 1{q^2}\overline{\Gamma }_{12} \\ 
-\frac 1{q^2}\overline{\Gamma }_{21} & \frac 1{q^2}\left( \overline{\Gamma }%
_{11}+(q^2-1)\overline{\Gamma }_{22}\right)
\end{array}
\right) }{\det_q(\overline{\Gamma })},  \label{gbinv} \\
\Gamma ^{-1} &=&\frac{\left( 
\begin{array}{cc}
q^2\Gamma _{22}-(q^2-1)\Gamma _{11} & -q^2\Gamma _{12} \\ 
-q^2\Gamma _{21} & \Gamma _{11}
\end{array}
\right) }{\det_{\frac 1q}(\Gamma ^T)}  \label{ginv}
\end{eqnarray}
and (\ref{gg}).

By expressing $q$ in terms of two parameters $\lambda $ and $\hbar $
according to $q=e^{\hbar \lambda }$ we can define two distinct limits: $%
\hbar \rightarrow 0$ which we refer to as the ''classical limit'' and $%
\lambda \rightarrow 0$ which we call ''canonical limit''. It is the \textbf{%
latter }and not the former that leads to Poincar\'e algebra. For this to
happen we will need to make $\Gamma $ depend on $\lambda $ in a certain way
as we explain in Sec.\ref{canlim}.

The outline of this paper is as follows. In Sec.\ref{claslim}, \ref{canlim}
we examine the limits $\hbar \rightarrow 0$ and $\lambda \rightarrow 0$,
respectively. The covariance of the algebra (\ref{pp}-\ref{pg}) under $%
SL_q(2,C)$ transformations is shown in Sec. \ref{covar}. The analogue of the
Pauli-Lubanski vector is discussed in Sec. \ref{PL}. We use it along with $P$
to construct two Casimir operators in Sec. \ref{casimirs}. These Casimir
operators are invariant under the action of the quantum Lorentz group. In
Sec. \ref{set} we find a set of commuting operators which can be used to
construct the representation of our algebra. Some concluding remarks are
given in Sec. \ref{concl}.

\section{Classical limit\label{claslim}}

We begin with the classical limit. Since the only dependence on $\hbar $
comes from $R-$ matrix, we need to know its classical limit\footnote{%
Althougth it looks like $R=e^{-i\hbar r},$ it is not so. Starting from the
third order in $\hbar $ this equality is violated but only in $3-2$ entry of
the matrix.} 
\begin{equation}
\underset{12}{R}\;\underset{\hbar \rightarrow 0}{\rightarrow }\Bbb{I+(}%
-i\hbar )r+\frac{(-i\hbar )^2}2r^2+O(\hbar ^3),  \label{Rser}
\end{equation}
where $r$ is the classical $r-$ matrix 
\begin{equation}
r=\frac{i\lambda }2\left( 
\begin{array}{cccc}
1 &  &  &  \\ 
& -1 &  &  \\ 
& 4 & -1 &  \\ 
&  &  & 1
\end{array}
\right)   \label{r}
\end{equation}
which satisfies the classical Yang-Baxter equation. Then the algebra (\ref
{pp}-\ref{pg}) becomes\footnote{%
From now on we use small letters for classical variables and capital letters
for quantum operators.} 
\begin{eqnarray}
\left\{ \underset{1}{p},\underset{2}{p}\right\}  &=&\left( r\underset{1}{p}%
\underset{2}{p}+\underset{1}{p}\underset{2}{p}r^{\dagger }-\underset{2}{p}%
r^{\dagger }\underset{1}{p}-\underset{1}{p}r\underset{2}{p}\right) ,
\label{ppc} \\
\left\{ \underset{1}{\gamma },\underset{2}{\gamma }\right\}  &=&\left(
r^{\dagger }\underset{1}{\gamma }\underset{2}{\gamma }+\underset{1}{\gamma }%
\underset{2}{\gamma }r-\underset{2}{\gamma }r\underset{1}{\gamma }-\underset{%
1}{\gamma }r^{\dagger }\underset{2}{\gamma }\right) ,  \label{ggc} \\
\left\{ \underset{1}{\gamma },\overline{\underset{2}{\gamma }}\right\} 
&=&\left( r\underset{1}{\gamma }\overline{\underset{2}{\gamma }}+\underset{1%
}{\gamma }\overline{\underset{2}{\gamma }}r-\overline{\underset{2}{\gamma }}r%
\underset{1}{\gamma }-\underset{1}{\gamma }r\overline{\underset{2}{\gamma }}%
\right) ,  \label{ggbc} \\
\left\{ \underset{1}{p},\underset{2}{\gamma }\right\}  &=&\left( r^{\dagger }%
\underset{1}{p}\underset{2}{\gamma }+\underset{1}{p}\underset{2}{\gamma }r-%
\underset{2}{\gamma }r^{\dagger }\underset{1}{p}-\underset{1}{p}r^{\dagger }%
\underset{2}{\gamma }\right) ,  \label{pgc}
\end{eqnarray}
where $\left\{ ,\right\} $ denotes a Poisson bracket which is obtained from
a commutator as follows: 
\begin{equation*}
\frac{\left[ ,\right] }{i\hbar }\underset{\hbar \rightarrow 0}{\rightarrow }%
\left\{ ,\right\} .
\end{equation*}
Obviously, this is not Poincar\'e algebra. It is instead the Poisson bracket
algebra discussed in \cite{ssy95}, \cite{sy95}. It was shown to be
consistent with Jacobi identities, hermitian conjugation, as well as being
Lie-Poisson with respect to Lorentz transformations.

\section{Canonical limit\label{canlim}}

We next discuss the canonical limit. Now we write $\Gamma =e^{i\lambda J}$
and $\overline{\Gamma }=e^{i\lambda J^{\dagger }}.$ Then in the limit $%
\lambda \rightarrow 0:$%
\begin{equation*}
\left[ \underset{1}{P},\underset{2}{P}\right] =0,\quad \left[ \underset{1}{J}%
,\underset{2}{J}^{\dagger }\right] =0,\quad \left[ \underset{1}{J},\underset{%
2}{J}\right] =2i\hbar \Pi \left( \underset{2}{J}-\underset{1}{J}\right)
,\quad \left[ \underset{1}{P},\underset{2}{J}\right] =i\hbar \underset{1}{P}%
\left( 2\Pi -\Bbb{I}\right) ,
\end{equation*}
where 
\begin{equation*}
\Pi =\left( 
\begin{array}{llll}
1 &  &  &  \\ 
&  & 1 &  \\ 
& 1 &  &  \\ 
&  &  & 1
\end{array}
\right)
\end{equation*}
is the permutation operator $\underset{1}{J}$ $\Pi =\Pi $ $\underset{2}{J},$ 
$\Pi $ $\underset{1}{J}=\underset{2}{J}$ $\Pi .$ The above algebra is the
Poincar\'e algebra. To see this we need only express 
\begin{eqnarray}
P &=&\left( 
\begin{array}{cc}
P_{11} & P_{12} \\ 
P_{21} & P_{22}
\end{array}
\right) =-\Bbb{I}P_0+\sigma _kP_k=\left( 
\begin{array}{cc}
-P_0+P_3 & P_1-i~P_2 \\ 
P_1+i~P_2 & -P_0-P_3
\end{array}
\right) ,  \label{P} \\
J &=&\sigma _k\left( J_{k0}-\frac i2\epsilon _{kmn}J_{mn}\right) =  \label{J}
\\
&=&\left( 
\begin{array}{cc}
-iJ_{12}+J_{30} & -iJ_{23}-iJ_{20}-J_{31}+J_{10} \\ 
-iJ_{23}+iJ_{20}+J_{31}+J_{10} & iJ_{12}-J_{30}
\end{array}
\right) ,  \notag
\end{eqnarray}
$\sigma _k$ being Pauli matrices, and then we arrive at the usual form: 
\begin{eqnarray*}
\left[ P_\mu ,P_\nu \right] &=&0,\quad \left[ P_\mu ,J_{\nu \rho }\right]
=i\hbar \left( \eta _{\mu \rho }P_\nu -\eta _{\mu \nu }P_\rho \right) , \\
\left[ J_{\mu \nu },J_{\rho \sigma }\right] &=&i\hbar \left( \eta _{\mu \rho
}J_{\nu \sigma }+\eta _{\nu \sigma }J_{\mu \rho }+\eta _{\mu \sigma }J_{\rho
\nu }+\eta _{\nu \rho }J_{\sigma \mu }\right) , \\
\eta &=&diag\left( -1,1,1,1\right) .
\end{eqnarray*}
Thus the algebra $\left( \text{\ref{pp}-\ref{pg}}\right) $ is a
two-parameter deformation of the Poincar\'e algebra.

\section{$SL_q(2,C)$ covariance\label{covar}}

It remains to show that $\left( \text{\ref{pp}-\ref{pg}}\right) $ is
covariant with respect to $SL_q(2,C)$ transformations. A $SL_q(2,C)$ matrix $%
T$ is defined using the commutational relations\footnote{%
In the canonical limit the matrix elements of $T$, of course, commute. In
the classical limit, however, we get non-trivial Poisson brackets for the
classical variables $t$ associated with $T$: 
\begin{equation*}
\left\{ \underset{1}{t},\underset{2}{t}\right\} =\left[ r,\underset{1}{t}%
\underset{2}{t}\right]  \label{sggc}
\end{equation*}
\begin{equation*}
\left\{ \underset{1}{t},\overline{\underset{2}{t}}\right\} =\left[ r,%
\underset{1}{t}\overline{\underset{2}{t}}\right]  \label{sggbc}
\end{equation*}
\begin{equation*}
\left\{ \overline{\underset{1}{t}},\overline{\underset{2}{t}}\right\}
=\left[ r,\overline{\underset{1}{t}}\,\overline{\underset{2}{t}}\right]
\label{sgbgbc}
\end{equation*}
here $t,$ $\overline{t}$ $\in SL(2,C)$, $\overline{t}=t^{\dagger -1}$}: 
\begin{equation}
\underset{12}{R}\underset{1}{T}\underset{2}{T}=\underset{2}{T}\underset{1}{T}%
\underset{12}{R}  \label{tt}
\end{equation}
which are well known (for review see e.g. \cite{takh90}). The action of $%
SL_q(2,C)$ on the operators $P$ and $\Gamma $ involves $T$ as well as its
hermitian conjugate $T^{\dagger }=\overline{T}^{-1}.$ We must then specify
the commutational relations involving $\overline{T}:$ 
\begin{eqnarray}
\underset{12}{R}\underset{1}{T}\underset{2}{\overline{T}} &=&\underset{2}{%
\overline{T}}\underset{1}{T}\underset{12}{R},  \label{ttb} \\
\underset{12}{R}~\underset{1}{\overline{T}}\underset{2}{~\overline{T}} &=&%
\underset{2}{\overline{T}}\underset{1}{~\overline{T}}~\underset{12}{R}.
\label{tbtb}
\end{eqnarray}
Since the ''deformed'' determinants of $T$ and $\overline{T}$ commute with
everything we can set them equal to one: 
\begin{eqnarray}
\det {}_{\frac 1{\sqrt{q}}}(T) &=&T_{11}T_{22}-qT_{12}T_{21}=1,  \label{dett}
\\
\det {}_{\frac 1{\sqrt{q}}}(\overline{T}) &=&\overline{T}_{11}\overline{T}%
_{22}-q\overline{T}_{12}\overline{T}_{21}=1.  \label{dettb}
\end{eqnarray}
Just like (\ref{c3}-\ref{c4}), (\ref{dett}-\ref{dettb}) are not independent
and one can be obtained from the other by hermitian conjugation with the
help of 
\begin{eqnarray}
T^{-1} &=&\frac{\left( 
\begin{array}{cc}
T_{22} & -\frac 1qT_{12} \\ 
-qT_{21} & T_{11}
\end{array}
\right) }{\det {}_{\frac 1{\sqrt{q}}}(T)},  \label{tinv} \\
\overline{T}^{-1} &=&\frac{\left( 
\begin{array}{cc}
\overline{T}_{22} & -\frac 1q\overline{T}_{12} \\ 
-q\overline{T}_{21} & \overline{T}_{11}
\end{array}
\right) }{\det {}_{\frac 1{\sqrt{q}}}(\overline{T})}.  \label{tbinv}
\end{eqnarray}
Under the action of $SL_q(2,C)$ the operators $P,$ $\Gamma $ and $\overline{%
\Gamma }$ transform according to 
\begin{eqnarray}
P &\rightarrow &P^{\prime }=\overline{T}PT^{-1},  \label{sym} \\
\Gamma &\rightarrow &\Gamma ^{\prime }=T\Gamma T^{-1},  \notag \\
\overline{\Gamma } &\rightarrow &\overline{\Gamma }^{\prime }=\overline{T}~%
\overline{\Gamma }~\overline{T}^{-1}.  \notag
\end{eqnarray}
The first relation states that $P$ transforms as a vector, while the latter
two define the adjoint action for $SL_q(2,C).$ Assuming, as usual, that the
the matrix elements of $\left( P,\Gamma ,\overline{\Gamma }\right) $ commute
with those of $\left( T,\overline{T}\right) $ it is not hard to show that $%
\left( \text{\ref{pp}-\ref{pg}}\right) $ are covariant under $\left( \text{%
\ref{sym}}\right) .$ For example, 
\begin{eqnarray*}
\underset{12}{R}\underset{1}{P}\underset{12}{R}^{-1}\underset{2}{P}
&\rightarrow &\underset{12}{R}\underset{1}{P}^{\prime }\underset{12}{R}^{-1}%
\underset{2}{P}^{\prime }=\underset{2}{\overline{T}}\,\underset{1}{\overline{%
T}}\underset{12}{R}\underset{1}{P}\underset{12}{R}^{-1}\underset{2}{P}%
\underset{1}{T}^{-1}\underset{2}{T}^{-1}= \\
\ &&\underset{2}{\overline{T}}\,\underset{1}{\overline{T}}\underset{2}{P}%
\underset{21}{R}^{-1}\underset{1}{P}\underset{21}{R}\underset{1}{T}^{-1}%
\underset{2}{T}^{-1}=\underset{2}{P}^{\prime }\underset{21}{R}^{-1}\underset{%
1}{P}^{\prime }\underset{21}{R}.
\end{eqnarray*}

\section{Pauli-Lubanski vector\label{PL}}

We obtained the following classical deformation of Pauli-Lubanski vector in 
\cite{ssy95}, \cite{sy95}: 
\begin{equation}
w=\frac 1{2\lambda }\left( \overline{\gamma }^{-1}p\gamma -p\right) =\left( 
\begin{array}{cc}
w_{11} & w_{12} \\ 
w_{21} & w_{22}
\end{array}
\right) =\left( 
\begin{array}{cc}
-w_0+w_3 & w_1-i~w_2 \\ 
w_1+i~w_2 & -w_0-w_3
\end{array}
\right)  \label{wc}
\end{equation}
which transforms as a Lorentz vector and reduces to the standard
Pauli-Lubanski vector when $\lambda \rightarrow 0.$ The Poisson brackets for 
$w$ are: 
\begin{eqnarray}
\left\{ \underset{1}{w},\underset{2}{\gamma }\right\} &=&r^{\dagger }%
\underset{1}{w}\underset{2}{\gamma }+\underset{1}{w}\underset{2}{\gamma }r-%
\underset{2}{\gamma }r^{\dagger }\underset{1}{w}-\underset{1}{w}r^{\dagger }%
\underset{2}{\gamma },  \label{wgc} \\
\left\{ \underset{1}{w},\underset{2}{\overline{\gamma }}\right\} &=&r%
\underset{1}{w}\underset{2}{\overline{\gamma }}+\underset{1}{w}\underset{2}{%
\overline{\gamma }}r-\underset{2}{\overline{\gamma }}r^{\dagger }\underset{1%
}{w}-\underset{1}{w}r\underset{2}{\overline{\gamma }},  \label{wgbc} \\
\left\{ \underset{1}{w},\underset{2}{p}\right\} &=&r\underset{1}{w}\underset{%
2}{p}+\underset{1}{w}\underset{2}{p}r^{\dagger }-\underset{2}{p}r^{\dagger }%
\underset{1}{w}-\underset{1}{w}r\underset{2}{p},  \label{wpc} \\
\left\{ \underset{1}{w},\underset{2}{w}\right\} &=&r\underset{1}{w}\underset{%
2}{w}+\underset{1}{w}\underset{2}{w}r^{\dagger }-\underset{2}{w}r^{\dagger }%
\underset{1}{w}-\underset{1}{w}r\underset{2}{w}-i\Pi \left( \underset{1}{w}%
\underset{2}{p}-\underset{2}{w}\underset{1}{p}\right) .  \label{wwc}
\end{eqnarray}
Because of the last term in the last Poisson bracket the quantum algebra for 
$w$ will depend on $\lambda $ not just through $q=e^{\hbar \lambda }$ but
also separately.

Let us define the quantum analog of $w$ according to: 
\begin{equation}
W=a\left( \beta \overline{\Gamma }^{-1}P\Gamma -P\right) ,  \label{w}
\end{equation}
where $W^{\dagger }=W,\quad a\underset{\hbar \rightarrow 0}{\rightarrow }%
\frac 1{2\lambda },$ $\beta \underset{\hbar \rightarrow 0}{\rightarrow }1.$
Using (\ref{w}) and (\ref{pp}-\ref{pg}) we get: 
\begin{eqnarray}
\underset{21}{R}^{-1}\underset{1}{W}\underset{21}{R}\underset{2}{\Gamma } &=&%
\underset{2}{\Gamma }\underset{21}{R}^{-1}\underset{1}{W}\underset{12}{R}%
^{-1},  \label{wg} \\
\underset{12}{R}\underset{1}{W}\underset{12}{R}^{-1}\underset{2}{\overline{%
\Gamma }} &=&\underset{2}{\overline{\Gamma }}\underset{21}{R}^{-1}\underset{1%
}{W}\underset{12}{R}^{-1},  \label{wgb} \\
\underset{12}{R}\underset{1}{W}\underset{12}{R}^{-1}\underset{2}{P} &=&%
\underset{2}{P}\underset{21}{R}^{-1}\underset{1}{W}\underset{21}{R},
\label{pw} \\
\underset{12}{R}(\underset{1}{W}+a\underset{1}{P})\underset{12}{R}^{-1}%
\underset{2}{W} &=&\underset{2}{W}\underset{21}{R}^{-1}(\underset{1}{W}+a%
\underset{1}{P})~\underset{21}{R}.  \label{ww}
\end{eqnarray}
$\left( \text{\ref{wg}-\ref{ww}}\right) $ are covariant with respect to $%
\left( \text{\ref{sym}}\right) $ and $W$ transforms as an $SL_q(2,C)$
vector: 
\begin{equation}
W\rightarrow W^{\prime }=\overline{T}WT^{-1}.  \label{wsym}
\end{equation}
Notice that in order to have correct classical limit $\hbar \rightarrow 0$
we have to include not only $q,$ but also $\lambda $ (through $a$) in the
commutational relations. Also note that $\left( \text{\ref{wg}-\ref{ww}}%
\right) $ are satisfied identically for any $\beta $ from $\left( \text{\ref
{w}}\right) $.

To recover the usual Pauli-Lubanski vector we assume that the canonical
limit of $a$ and $\beta $ is identical to the classical limit, i.e. $a%
\underset{\lambda \rightarrow 0}{\rightarrow }\frac 1{2\lambda },$ $\beta 
\underset{\lambda \rightarrow 0}{\rightarrow }1.$ Then using $\Gamma
=e^{i\lambda J}$ and $\overline{\Gamma }=e^{i\lambda J^{\dagger }}$ we get $%
W=-\Bbb{I}W_0+\sigma _kW_k,$ where 
\begin{equation*}
W_0=-\frac 12\epsilon _{kmn}P_kJ_{mn},\quad W_k=-\frac 12\epsilon
_{kmn}P_0J_{mn}-\epsilon _{lmk}P_lJ_{m0}.
\end{equation*}
In 4-vector notations it can be rewritten in the form: 
\begin{equation*}
W_\beta =-\frac 12\epsilon _\beta ^{~\mu \nu \rho }J_{\mu \nu }P_\rho .
\end{equation*}

\section{Casimirs\label{casimirs}}

Like with the Poincar\`e algebra, we can construct two Casimir operators.
Both are invariant under $SL_q(2,C)$ transformations as we shall make
evident. The following quadratic combination is an analogue of $P^2$ and is
one such Casimir: 
\begin{equation}
\mathcal{C}_1=-det_q(P)=-<P,P>_q=(P,P)_q=-\ P_0^2+\frac{P_1^2}{q^2}+\frac{%
P_2^2}{q^2}+~P_3^2\underset{q\rightarrow 1}{\rightarrow }P^2,  \label{c1}
\end{equation}
where 
\begin{equation}
<A,B>_q=A_{11}B_{22}-\frac 1{q^2}A_{12}B_{21}  \label{<>}
\end{equation}
is a combination of matrix elements which often occurs in this algebra so we
reserve a special notation for it. 
\begin{equation}
(A,B)_q=-\frac 1{q^2+1}Tr_q\left( A\widetilde{B}\right)   \label{()}
\end{equation}
is a deformed scalar product of 4-vectors $A$ and $B$ represented as $%
2\times 2$ hermitian matrices. $\widetilde{B}$ is a deformed adjugate of $B-$
matrix. By this we mean that $\widetilde{B}$ has the properties: 
\begin{eqnarray*}
B\widetilde{B} &=&\widetilde{B}B\sim \Bbb{I},\quad  \\
\widetilde{B}\underset{q\rightarrow 1}{\rightarrow }\widetilde{b} &=&\left( 
\begin{array}{ll}
b_{22} & -b_{12} \\ 
-b_{21} & b_{11}
\end{array}
\right) ,\quad \widetilde{b}b=b\widetilde{b}=\det (b)\Bbb{I}.
\end{eqnarray*}
$Tr_q$ is a deformed trace defined by
\begin{equation*}
Tr_q(A)=A_{11}+q^2A_{22}.
\end{equation*}
The adjugate of $P$ has the following form: 
\begin{equation}
\widetilde{P}=\left( 
\begin{array}{cc}
P_{22} & -\frac 1{q^2}P_{12} \\ 
-\frac 1{q^2}P_{21} & \frac 1{q^2}\left( P_{11}+(q^2-1)P_{22}\right) 
\end{array}
\right) .  \label{Pad}
\end{equation}
It was proved directly using Mathematica that $\widetilde{P}$ has the same
transformation properties\footnote{%
When $\det_q(P)\neq 0$ 
\begin{equation*}
\widetilde{P}=det_q(P)P^{-1}.
\end{equation*}
Notice though that $\widetilde{P}$ exists even when $P^{-1}$ does not, that
is for 0-length 4-vectors. To find $\widetilde{P}$ we used an anzatz where
each element of the matrix is a general linear combination of all $P_{ij}.$ $%
\widetilde{P}$ is uniquely determined up to a factor which goes to 1 when $%
q\rightarrow 1.$} as $\overline{P}:$%
\begin{equation*}
\widetilde{P}^{\prime }\rightarrow T\widetilde{P}\overline{T}^{-1}.
\end{equation*}

We found another Casimir corresponding to the square of the Pauli-Lubanski
vector. It can be expressed in different ways: 
\begin{eqnarray}
\mathcal{C}_{2a} &=&a\left( <P,W>_q+q^2<W,P>_q\right) =  \label{c2a} \\
\  &=&aTr_q\left( W\widetilde{P}\right) =aTr_q\left( \widetilde{P}W\right)
=-a(q^2+1)(W,P)_q,  \notag  \label{c2b} \\
\mathcal{C}_{2b} &=&aTr_q\left( P\widetilde{W}\right) =aTr_q\left( 
\widetilde{W}P\right) =-a(q^2+1)(P,W)_q,  \label{c2b} \\
\mathcal{C}_{2c} &=&-<W,W-a(q^2-1)P>_q=(W,W)_q,  \label{c2c}
\end{eqnarray}
Here\footnote{%
Again, when $<W,W-a(q^2-1)P>_q\neq 0$%
\begin{equation*}
\widetilde{W}=<W,W-a(q^2-1)P>_qW^{-1}.
\end{equation*}
To find $\widetilde{W}$ we used an anzatz where each element of the matrix
is a general linear combination of all $W_{ij}$ and $P_{ij}.$ It was
uniquely determined up to a factor which goes to 1 when $q\rightarrow 1.$} 
\begin{equation}
\widetilde{W}=\left( 
\begin{array}{cc}
W_{22} & -\frac{W_{12}}{q^2} \\ 
-\frac{W_{21}}{q^2} & \frac{W_{11}+(q^2-1)W_{22}}{q^2}
\end{array}
\right) -a(q^2-1)\widetilde{P}  \label{Wad}
\end{equation}
is an adjugate matrix of $W$ and direct calculations show that $\widetilde{W}
$ transforms like $\overline{W}:$%
\begin{equation*}
\widetilde{W}^{\prime }\rightarrow T\widetilde{W}\overline{T}^{-1}.
\end{equation*}
One can check using (\ref{w}) that $\mathcal{C}_{2a},$ $\mathcal{C}_{2b}$
and $\mathcal{C}_{2c}$ are in fact equal to each other up to $\mathcal{C}_1$%
: 
\begin{equation*}
\mathcal{C}_{2a}=\mathcal{C}_{2b}-a^2(q^4-1)\mathcal{C}_1=\mathcal{C}_{2c}+%
\frac{a^2(q^6-\beta ^2)}{q^6}\mathcal{C}_1.
\end{equation*}
$\mathcal{C}_{2a}\ $and $\mathcal{C}_{2b}$ have the same classical limit,
namely $\emph{c}_{2a}=-\frac{(p,w)}\lambda ,$ where $\left( ,\right) =\left(
,\right) _{q=1}.$ The classical limit of $\mathcal{C}_{2c}$ is $\emph{c}%
_{2c}=w^2.$ At first glance one might think that $\emph{c}_{2a}\neq \emph{c}%
_{2c}$ but it is easy to prove that $\emph{c}_{2a}$ and $\emph{c}_{2c}$ are
really equal to each other. For this we write 
\begin{equation*}
w=\frac 1{2\lambda }(b-p),
\end{equation*}
where $b=\overline{\gamma }^{-1}p\gamma $. Then $\det (b)=\det (p)$ because $%
\det (\gamma )=\det (\overline{\gamma })=1.$ Furthermore 
\begin{equation*}
2\lambda (p,w)=(p,b)+\det (p)
\end{equation*}
using $(p,p)=-\det (p),$ and thus 
\begin{equation*}
w^2=\frac 1{4\lambda ^2}(-\det (b)-\det (p)-(b,p)-(p,b))=\frac 1{2\lambda
^2}(-\det (p)-(p,b))=-\frac{(p,w)}\lambda .
\end{equation*}

Let us show that our deformed scalar product is really invariant and hence
that $\mathcal{C}_1$ and $\mathcal{C}_2$ are invariant. Suppose that under $%
SL_q(2,C)$ transformations 
\begin{equation*}
A^{\prime }\rightarrow \overline{T}AT^{-1},\quad \widetilde{B}^{\prime
}\rightarrow T\widetilde{B}\overline{T}^{-1}
\end{equation*}
and matrix elements of $T$ commute with $A$ and $\widetilde{B}.$ Then 
\begin{eqnarray*}
Tr_q\left( A^{\prime }\widetilde{B}^{\prime }\right)  &=&Tr_q\left( 
\overline{T}A\widetilde{B}\overline{T}^{-1}\right) =\left( \overline{T}_{1k}%
\overline{T}_{n1}^{-1}+q^2\overline{T}_{2k}\overline{T}_{n2}^{-1}\right)
A_{km}\widetilde{B}_{mn}= \\
\  &=&\delta _{kn}^qA_{km}\widetilde{B}_{mn}=Tr_q\left( A\widetilde{B}%
\right) ,\quad where\quad \delta ^q=\left( 
\begin{array}{ll}
1 & 0 \\ 
0 & q^2
\end{array}
\right) .
\end{eqnarray*}
Notice that in all cases discussed above $Tr_q\left( A\widetilde{B}\right)
=Tr_q\left( \widetilde{B}A\right) $ which is rather unexpected for such a
deformed algebra. On the other hand, $Tr_q\left( A\widetilde{B}\right) \neq
Tr_q\left( B\widetilde{A}\right) ,$ that is the scalar product is not
symmetric  in general.

\section{Complete set of mutually commuting operators\label{set}}

In addition to the Casimirs the following operators can be included into a
complete set of mutually commuting operators necessary to construct a
representation: 
\begin{eqnarray}
\mathcal{K}_1 &=&Tr_q(P)=P_{11}+q^2P_{22},  \label{c5} \\
\mathcal{K}_2 &=&Tr_q(W)=W_{11}+q^2W_{22}=a\left( qTr_q(P\Omega )-\mathcal{K}%
_1\right) ,  \label{c6} \\
\mathcal{K}_3 &=&P_{11}\quad or\quad P_{22},  \label{c7} \\
\mathcal{K}_4 &=&<\Gamma ,\overline{\Gamma }>_q=\Omega _{11},  \label{k4}
\end{eqnarray}
where 
\begin{equation}
\Omega =\Gamma \overline{\Gamma }^{-1}=\Gamma \Gamma ^{\dagger }=\Omega
^{\dagger }.  \label{sigma}
\end{equation}
When constructing a representation, it is more convenient to choose $%
\mathcal{K}_3=P_{22}$ because it has simpler commutational relations with
other operators.

$\Omega $ has very peculiar property: it has the same commutational
relations with all matrices $Z=P,$ $\Gamma ,$ $\overline{\Gamma },$ $W$ or $%
\Omega $: 
\begin{equation}
\underset{1}{Z}\underset{21}{R}\underset{2}{\Omega }\underset{12}{R}=%
\underset{21}{R}\underset{2}{\Omega }\underset{12}{R}\underset{1}{Z}.
\label{sigmacr}
\end{equation}
It follows from here, by the way, that components of $\Omega $ generate a
universal enveloping algebra.

In canonical limit $\mathcal{K}_1$ and $\mathcal{K}_2$ go to $-2P_0$ and $%
-2W_0,$ respectively, while 
\begin{equation}
\frac{\Omega -\Bbb{I}}\lambda \underset{\lambda \rightarrow 0}{\rightarrow }%
\left( 
\begin{array}{ll}
J_{12} & J_{23}-iJ_{31} \\ 
J_{23}+iJ_{31} & -J_{12}
\end{array}
\right)  \label{omegalim}
\end{equation}
is a traceless hermitian matrix corresponding to space components of angular
momentum. Therefore $\mathcal{K}_4$ is associated with the component of the
angular momentum along the third axis $J_{12}.$

It was checked that no other linear or quadratic combination of $P,$ $W,$ $%
\Gamma $ and $\overline{\Gamma }$ can be included in this set.

Another possible set of mutually commuting operators can be obtained if we 
\textit{replace} $\mathcal{K}_3$ by 
\begin{equation}
\mathcal{K}_5=Tr\left( \Omega \right) \underset{\lambda \rightarrow 0}{%
\rightarrow }2\left( 1+2\lambda ^2\left( J_{12}^2+J_{23}^2+J_{31}^2\right)
\right) +o\left( \lambda ^2\right)  \label{k5}
\end{equation}
which corresponds to the square of 3-vector of angular momentum (\textit{%
replace}, because no linear combination of $P_{11}$ and $P_{22},$ except for 
$\mathcal{\QTR{mathcal}{K}}_1,$ commutes with $\mathcal{K}_5$ ). Therefore
in the set we can have either an analogue of $P_3$ or $\overrightarrow{J}^2$
but not both\footnote{%
At first glance one might think that $\det (\Omega )$ can also be included
in either of two sets but it turned out that it is not independent: 
\begin{equation*}
\det (\Omega )=\frac 1{q^2}\left( 1+\left( q^2-1\right) \mathcal{K}%
_4^2\right) .
\end{equation*}
}.

For the sake of space we do not write out all commutational relations
between $\mathcal{K}_i$ and elements of the algebra but only some simple
ones: 
\begin{eqnarray*}
\lbrack \mathcal{\QTR{mathcal}{K}}_1,P_{ij}] &=&0,\quad [\mathcal{%
\QTR{mathcal}{K}}_2,P_{ij}]=0, \\
\lbrack \mathcal{\QTR{mathcal}{K}}_4,Z_{11}] &=&[\mathcal{\QTR{mathcal}{K}}%
_4,Z_{22}]=0,\quad \mathcal{\QTR{mathcal}{K}}_4Z_{12}=\frac 1{q^2}Z_{12}%
\mathcal{\QTR{mathcal}{K}}_4,\quad \mathcal{\QTR{mathcal}{K}}%
_4Z_{21}=q^2Z_{21}\mathcal{\QTR{mathcal}{K}}_4,
\end{eqnarray*}
where $Z$ can be: $P,$ $W,$ $\Gamma $ , $\overline{\Gamma }$ or $\Omega $.
We shall give a complete set of commutational relations along with
representations of the algebra in a forthcoming article.

\section{Conclusion\label{concl}}

Here we remark on possible extension of our algebra. In \cite{ssy95} we
deformed the canonical Poisson brackets for a relativistic particle so that
the $SL(2,C)$ Poisson-Lie group with Poisson structure given in footnote$^3$
had a Poisson action on the classical observables. The most non-trivial part
of \cite{ssy95} was to deform the canonical Poisson brackets between momenta
and coordinates $\left\{ x_\mu ,p_\nu \right\} =\eta _{\mu \nu }.$ We
obtained 
\begin{equation}
\left\{ \underset{1}{x},\underset{2}{p}\right\} =r\underset{1}{x}\underset{2%
}{p}+\underset{1}{x}\underset{2}{p}r^{\dagger }-\underset{2}{p}r\underset{1}{%
x}-\underset{1}{x}r^{\dagger }\underset{2}{p}-\Pi (\underset{1}{f^{\dagger }}%
+\underset{2}{f})  \label{xpc}
\end{equation}
where $f=\exp (i\sin ^{-1}\left( \lambda xp\right) ).$ (\ref{xpc}) was shown
to be covariant with respect to $x\rightarrow tx\overline{t}^{-1}$, $%
p\rightarrow \overline{t}pt^{-1}.$ It would be interesting to extend quantum
algebra (\ref{pp}-\ref{pg}) to include the space-time operator $X$ . Upon
choosing the Hamiltonian to be equal to $\mathcal{C}_1$we could then write
down the corresponding Klein-Gordon equation giving the dynamics of a
''free'' particle, the classical analogue of which was discussed in \cite
{ssy95}, \cite{sy95}.

\section{Acknowledgments}

This work was in part supported by Department of Energy, USA, under contract
number DE-FG0584ER40141. I.Y. was also supported by Graduate Council
Research Fellowship from the University of Alabama.

\bibliographystyle{unsrt}
\bibliography{btxdoc,qg1}

\end{document}